\documentclass[12pt,reqno]{article}
\usepackage[]{fontenc}
\usepackage[latin1]{inputenc}
\usepackage[francais]{babel}
\usepackage{a4,amsmath,amssymb,amsfonts,amsthm}
\usepackage{amsrefs}

\begin{document}
\bibliographystyle{plain}

% theorems etc.

 \newtheorem{theor}{Théorème}[section]
 \newtheorem{prop}[theor]{Proposition}
 \newtheorem{cor}[theor]{Corollaire}
 \newtheorem{lemme}[theor]{Lemme}
 \newtheorem{slemme}[theor]{sous-lemme}
 \newtheorem{defin}[theor]{Définition}
 \newtheorem{conj}[theor]{Conjecture}

% abbreviations

\gdef\Td{{\rm Td}}
\gdef\ch{{\rm ch}}
\gdef\CB{{\cal B}}
\def\CF{{\cal F}}
\gdef\CE{{\cal E}}
\gdef\CL{{\cal L}}
\gdef\CO{{\cal O}}
\gdef\CS{{\cal S}}
\gdef\CM{{\cal M}}
\gdef\CC{{\cal C}}
\gdef\CA{{\cal A}}
\gdef\CH{{\cal H}}
\gdef\CT{{\cal T}}
\def\M#1{\mathbb#1}     % fontes tableau noir (ensembles R, C, Z, ...)
\def\B#1{\bold#1}       % symboles math\'ematiques accentues
\def\C#1{\mathcal#1}    % lettres caligraphiees (mode math ou texte)
\def\E#1{\scr#1}        % autres fontes caligraphiees (mode math ou texte)
\def\mR{\M{R}}           
\def\mZ{\M{Z}}           
\def\mN{\M{N}}           
\def\mQ{\M{Q}}       
\def\mC{\M{C}}  
\def\Qb{{\overline{\mQ}}}
\def\Qmn{\mQ({\mu_n})}
\gdef\beginProof{\par{\bf Preuve. }}
\gdef\endProof{${\bf Q.E.D.}$\par}  
\gdef\ari#1{\widehat{#1}}
\gdef\mtr#1{\overline{#1}}
\gdef\c1{{{\rm c}_1}}
\gdef\ac1{\ari{\rm c}_1}
\gdef\aceq1{\ari{\rm c}_{{\rm eq},1}}
\gdef\acheq{\ari{\rm ch}_{{\rm eq}}}
\gdef\ach{\ari{\rm ch}}
\gdef\ra{\rightarrow}
\gdef\Spec{\rm Spec}
\gdef\Hom{{\rm Hom}}
\gdef\tbcov{{{\cal T}_{\rm BCOV}}}
\gdef\Vol{{\rm Vol}}
\gdef\NT{{\rm NT}}
\gdef\refeq#1{(\ref{#1})}
\def\Gr{{\rm Gr}}
\def\Gal{{\rm Gal}}
\def\T{{\rm T}}
\def\H{{\rm H}}
 \def\R{{\rm R}}
 \def\P{{\rm P}}
 \def\Pet{{\rm Pet}}
 \def\Tw{{\rm Tw}}
 \def\Jac{{\rm Jac}}
 \def\Cl{{\rm Cl}}
 \def\HDL{{\rm H}_{\rm dlb}}
 \def\reg{{\rm reg}}
 \def\mff{{\mathfrak f}}
 \def\mn{{\mu_n}}
\def\m2{{\mu_2}}
 \def\rk{{\rm rk}}
 \def\Zn{{\mZ/n\mZ}}
 \def\K{{\rm K}}
 \def\ul#1{\underline{#1}}
 \def\Tor{{\rm Tor}}
 \def\mfp{{\mathfrak p}}
\def\adeg{\ari{\rm deg}}
\def\op#1{\operatorname{#1}}
\def\Zar{{\rm Zar}}

% heading

\author{Vincent Maillot\footnote{CNRS, Centre de Math\'ematiques de Jussieu,  
     Universit\'e Paris 7 Denis Diderot,  
     Case Postale 7012,  
     2, place Jussieu,  
     F-75251 Paris Cedex 05, France}\ \  et Damian R\"ossler\footnote{CNRS, Centre de Math\'ematiques de Jussieu,  
     Universit\'e Paris 7 Denis Diderot,  
     Case Postale 7012,  
     2, place Jussieu,  
     F-75251 Paris Cedex 05, France}}
\title{Formes automorphes et théorèmes de Riemann-Roch arithmétiques}
\maketitle
\date
\abstract{Nous construisons trois familles de formes automorphes au moyen du théorème 
de Riemann-Roch arithmétique et de la formule de Lefschetz arithmétique. 
Deux de ces familles ont déjà été construites par Yoshikawa et notre construction 
met en lumière leur origine arithmétique.}
% beginning of the text

\parindent=0pt
\parskip=5pt

\section{Introduction}

Le but de ce texte est de donner une interprétation arithmétique et géométrique 
à trois familles de formes automorphes d'expression analytique.
Plus précisément, on démontre que ces formes 
automorphes sont algébriques et entières, lorsque les espaces-sous jacents 
ont des modèles entiers. 

La première est la famille de formes modulaires de Siegel construite par Yoshikawa dans 
\cite{Yoshikawa-Discriminant}.
Notre calcul démontre une version légèrement affaiblie d'une conjecture de Yoshikawa sur les 
coefficients de Fourier de ces formes modulaires. 
 La 
deuxième est la famille de formes modulaires d'Igusa \og produit des 
thêta constantes paires\fg.
La troisième est la famille de formes 
automorphes à coefficients sur certains espaces de modules de surfaces 
K3 à coefficients; lorsque l'involution est sans point fixe, elles 
coïncident avec certaines fonctions $\Phi$ de Borcherds 
(cf. \cite[Sec. 8]{Yoshikawa-K3-surfaces}). Cette famille de formes 
est construite par Yoshikawa dans 
\cite[Th. 5.2]{Yoshikawa-K3-surfaces}. Notre calcul démontre en particulier 
que les fonctions $\Phi$ ci-dessus sont d'origine arithmétique. 

Dans l'appendice, nous formulons une extension conjecturale de 
la formule de Lefschetz arithmétique, où des irrégularités sont 
autorisées sur les fibres finies. Cette formule n'est pas appliquée dans le présent 
texte mais elle représente un moyen théorique d'étudier la dégén\'erescence 
de la deuxième forme modulaire de Yoshikawa lorsqu'on considère 
une surface K3 avec involution définie sur un corps de nombres 
et ayant mauvaise réduction en certaines places finies.

Dans ce texte, nous utiliserons librement la terminologie et les résultats 
énoncés dans la section 4 de \cite{Koehler-Rossler-A-fixed-point} (article dans lequel la formule 
de Lefschetz arithmétique mentionnée plus haut est démontrée). 
Par ailleurs, nous utiliserons la terminologie et les résultats de 
\cite{Gillet-Soule-An-arithmetic} 
(article dans lequel le théorème de Riemann-Roch arithmétique en degré $1$ est 
démontré).

L'objet du présent texte est de présenter des calculs. Pour 
une introduction au théorème de Riemann-Roch arithmétique et 
à la formule de Lefschetz arithmétique, nous suggérons de 
consulter les articles originaux cités dans le dernier paragraphe, ainsi que 
\cite{Gillet-Rossler-Soule-An-arithmetic} ou encore 
les notes \cite{Rossler-Notes}.

{\bf Remerciements.} Une partie de ce travail a \'et\'e r\'ealis\'ee
alors que le premier auteur \'etait professeur invit\'e
au R.I.M.S. (Universit\'e de Kyoto). Il lui est tr\`es
agr\'eable de remercier cette institution pour son hospitalit\'e
et les conditions de travail exceptionnelles dont il a pu b\'en\'eficier. 
Nos remerciements vont \'egalement à K.-I. Yoshikawa, pour toutes les explications 
qu'il nous a fournies sur ses travaux.

 \section{Les formes modulaires de Yoshikawa de premier type}

Soit $S$ le spectre d'un anneau arithmétique. Soit 
$B$ une variété arithmétique sur $S$. Dans ce texte, on appellera 
{\it variété arithmétique} sur $S$ un schéma intègre 
et régulier, qui est quasi-projectif sur $S$. 
Soit $\pi:\CA\to B$ un schéma abélien de dimension relative $g$. 
Soit $h:\Theta\to B$ un morphisme lisse et propre de dimension relative
$g-1$ et $\theta:\Theta\hookrightarrow\CA$ une $B$-immersion 
fermée. On suppose que $\CO(\Theta)$ est un fibré 
relativement ample et que le degré de $\CO(\Theta)$ est 
$g!$ sur chaque fibre géométrique de $\CA/B$. Une hypothèse 
équivalente est que la caractéristique d'Euler de 
$\CO(\Theta)$ vaut $1$ sur chaque fibre géométrique de $\CA/B$. 

Nous noterons  $T\Theta:=\Omega^\vee_{\Theta/B}$ et 
$T\CA:=\Omega^\vee_{\CA/B}$. Nous écrirons $u:B\to\CA$ pour la section 
unité et $TA_0$ pour $u^*TA$. Nous écrirons aussi 
$\omega:=u^*\det(\Omega_{\CA/B})$. On note 
$\mtr{\CO}(\Theta)$ le fibré ${\CO}(\Theta)$ muni de sa métrique de 
Moret-Bailly (voir \cite[Par. 3.2]{Moret-Bailly-Sur}) et on pose  
$\mtr{L}:=\mtr{\CO}(\Theta)\otimes\pi^*u^*\mtr{\CO}(\Theta)^\vee$. La 
forme $2\pi\cdot \c1(\mtr{L})$ définit une structure de fibration Kählerienne sur 
$\CA$ au sens de \cite[Par. 1]{Bismut-Koehler-Higher}. 
% VERIFIER CE FAIT
Soit $N$ le fibré conormal de 
l'immersion $\theta$. 

Le {\it morphisme de Gauss} est défini de la manière suivante. 
Le morphisme naturel  $\T\Theta\hookrightarrow \theta^*\T\CA$ 
induit un morphisme $\Theta\to \Gr(g-1,\theta^*\T\CA)$. 
Utilisant les isomorphismes naturels $\Gr(g-1,\theta^*\T\CA)\simeq 
\theta^*\Gr(g-1,\T\CA)$, $\Gr(g-1,\pi^*\T\CA)\simeq 
\pi^*\Gr(g-1,\T\CA)$ et 
$\T A\simeq\pi^*\T\CA_0$, on obtient un morphisme naturel 
$\Theta\to h^*\Gr(g-1,\T\CA_0)$. Si l'on compose ce dernier avec la 
projection naturelle de $h^*\Gr(g-1,\T\CA_0)=\Gr(g-1,\T\CA_0)\times_B\Theta$ 
sur le premier facteur, on obtient le morphisme de Gauss 
$\gamma:\Theta\to\Gr(g-1,\T A_0)$. On note \makebox{$p:\P:=\Gr(g-1,\T A_0)\to B$} l'application structurale.  Pour la définition de $\Gr(\cdot,\cdot)$ voir 
\cite[App. B.5.7]{Fulton-Intersection}. On notera  
$$
{\cal E}: 0\to E\to p^*T\CA_0\to Q\to 0
$$
la suite exacte universelle sur $P$. Si l'on munit $p^*T\CA_0$ 
de la métrique image réciproque de celle de $T\CA_0$ et les fibrés $E$ et 
$Q$ des métriques induites, on obtient à partir de $\cal E$ une suite 
exacte métrisée que nous noterons $\mtr{\cal E}$. 

\begin{lemme}
Le morphisme de Gauss est génériquement fini de degré $g!$. 
\end{lemme}
\beginProof
Le fait que le morphisme de Gauss est génériquement fini 
(ou autrement dit, qu'il induit une extension finie de corps de fonctions 
$\kappa(\Theta)|\kappa(\Gr(g-1,\T A_0))$) est démontré dans 
\cite[Th. 4]{Abramovich-Subvarieties}. Pour calculer son degré, nous considérons 
 le calcul suivant dans la théorie de Chow de $\CA$:
 \begin{eqnarray*}
 g!&\stackrel{(1)}{=}&\pi_*(\c1(\CO(\Theta))^g)=
 \pi_*(\theta_*(1)\c1(\CO(\Theta))^{g-1})\stackrel{(2)}{=}
 h_*(\c1(\theta^*(\CO(\Theta)))^{g-1})\\
 &\stackrel{(3)}{=}&
 h_*(\c1(N^\vee)^{g-1})\stackrel{(4)}{=}
 p_*\gamma_*(\gamma^*(\c1(Q)^{g-1}))\stackrel{(5)}{=}
\deg(\gamma)p_*(\c1(Q)^{g-1})\\&=&\deg(\gamma).
\end{eqnarray*}
L'égalité (1) est justifiée par le théorème de Riemann-Roch 
(appliqué au morphisme $\pi$ et au fibré $\CO(\Theta)$), l'égalité 
(2) est justifiée par la formule de projection, l'égalité (3) est justifiée par 
la formule d'adjonction, l'égalité (4) est une conséquence la définition 
de $\P$ et pour finir (5) est une conséquence du fait que 
le degré de $\CO(1)$ sur un espace projectif au-dessus d'un corps 
vaut $1$.
\endProof

Nous appliquons maintenant le théorème de Riemann-Roch 
arithmétique à $\Theta$. 

\begin{lemme}
Les égalités suivantes  
\begin{eqnarray*}
\ac1(\R^\bullet h_*(\mtr{\CO}_\Theta))&=&
(-1)^g\ac1(\R^0 \pi_*(\mtr{\Omega}^g_\CA))+\log({g(g-2)\cdots (1+(1+(-1)^g)/2)\over 
(g-1)(g-3)\cdots(1+(1-(-1)^g)/2)})\\
&=&
(-1)^{g+1}\ac1(\mtr{\T A}_0)+\log({g(g-2)\cdots (1+(1+(-1)^g)/2)\over 
(g-1)(g-3)\cdots(1+(1-(-1)^g)/2)})
\end{eqnarray*}
sont vérifiées. 
\label{IDKLem}
\end{lemme}
\beginProof
Soit $M$ une variété Kählerienne de dimension $g$ et de forme 
de Kähler $\underline{\omega}$. Soit $k\geqslant 0$ et 
soit $\nu\in\H^k(M,\CO_M)$. 
On dispose de la formule suivante : 
\begin{equation}
 |\nu|^2_{L_2}:={ i^{k}(-1)^{k(k+1)/2}\over
 (2\pi)^g(g-k)!}\int_{M}\nu\wedge\mtr{\nu}\wedge 
 \ul{\omega}^{g-k}.
\label{metL2}
\end{equation}
Voir \cite[Par. 2.3]{Maillot-Rossler-On-the}. Par ailleurs, considérons la suite exacte longue de 
cohomologie
\begin{eqnarray*}
0&\to& \R^0 \pi_*\CO\stackrel{}{\to}\R^0 h_*\CO_\Theta\to 0\\
&\to&\R^1 \pi_*\CO
\stackrel{}{\to} \R^1 h_*\CO_\Theta\to 0\\
&\to&\dots\to\\
0&\to& \R^{g-1} \pi_*\CO
\stackrel{}{\to} \R^{g-1} h_*\CO_\Theta\to
\R^g \pi_*\CO(-\Theta)\stackrel{\star}{\longrightarrow}  \R^g \pi_*\CO\to 0
\end{eqnarray*}
née de la suite exacte
$$
0\to \CO(-\Theta)\to \CO\to \CO_\Theta\to 0. 
$$
On remarque aussi que la flèche $\star$ est un isomorphisme car 
$\R^g \pi_*\CO$ est localement libre de rang $1$. Toutes les flèches 
reliant deux objets non-nuls dans la suite exacte longue sont donc des isomorphismes. 
En particulier les faisceaux de cohomologie 
$\R^k h_*\CO_\Theta$ sont localement libres. 
De plus, en comparant la formule 
\refeq{metL2} sur les fibres de $\CA(\mC)$ et  sur les fibres 
de $\Theta(\mC)$, on conclut que pour tout entier $k$ tel que $0\leqslant k\leqslant g-1$,  on a 
\begin{equation}
\ac1(\R^k \pi_*(\mtr{\CO}_\CA))=\ac1(\R^k h_*(\mtr{\CO}_\Theta))-
\log(g-k).
\label{cpmet}
\end{equation}
On utilisé ici le fait que la forme de Kähler sur les fibres est donnée par 
$2\pi\cdot\c1(\mtr{\CO}(\Theta))$. 
Si l'on combine cette dernière égalité avec 
l'égalité 
\begin{eqnarray}
\ac1(\R^\bullet \pi_*(\mtr{\CO}_\CA))=0
\label{RRA}
\end{eqnarray}
on obtient la première égalité du lemme. L'égalité \refeq{RRA} est une conséquence 
immédiate du théorème de Riemann-Roch arithmétique 
appliqué à $\pi$ et $\mtr{\CO}_\CA$ et de l'annulation de la torsion 
analytique du fibré trivial d'une variété abélienne de dimension 
$\geqslant 2$ (cf. \cite[Par. 5, p. 173]{Ray-Singer-Analytic}). La deuxième égalité du lemme 
est une conséquence de la formule de projection.
\endProof

On suppose dorénavant que $\CO(\Theta)$ satisfait à l'hypothèse du 
dernier lemme. Le théorème de Riemann-Roch arithmétique 
appliqué à $h$ et au fibré hermitien trivial sur $\Theta$ 
donne l'égalité suivante dans $\ari{\rm CH}^{\leqslant 1}(B)_\mQ$:
\begin{eqnarray*}
\ach(R^\bullet h_*\mtr{\CO}_\Theta)-T(\Theta,\mtr{\CO}_\Theta)
&=&h_*(\ari{\Td}(\mtr{\T h}))-\int_{\Theta/B}\R(\T h)\Td(\T h)\\
&=& g!\ p_*(\ari{\Td}(\mtr{E}))- 
g!\ \int_{P/B}\R(E)\Td(E)\\
&=&
g!\ \int_{P/B}\Td^{-1}(\mtr{Q})\widetilde{\Td}(\mtr{\cal E})+g!\ p_*(\ari{\Td}^{-1}(\mtr{Q})
\ari{\Td}(p^*\mtr{\T\CA}_0))\\
&-& 
g!\int_{P/B}\R(E)\Td(E).
\end{eqnarray*}
% VERIFIER FORMULE D'ECHANGE
Remarquons à présent que pour tout $k\geqslant 0$, on a
$$
\ari{\Td}^{-1}(\mtr{Q})^{[k]}=-{(-1)^{k+1}\over (k+1)!}\ari{c}_1(\mtr{Q})^k
$$
(on rappelle que par définition $\Td^{-1}(x)=(1-\exp(-x))/x$). 
Par ailleurs, selon \cite[Par. 8.2]{Mourougane-Computations}
$$
p_*(\ari{c}_1(\mtr{Q})^{g}))=[\sum_{l=1}^{g-1}\CH_l]+\ari{c}_1(\mtr{\T\CA}_0);
$$
ici $\CH_l:=\sum_{k=1}^l{1\over k}$ est le $l$-ème nombre harmonique. 
On peut maintenant calculer 
\begin{eqnarray*}
&&g!\ p_*(\ari{\Td}^{-1}(\mtr{Q})
\ari{\Td}(p^*\mtr{\T\CA}_0))^{[\leqslant 1]}\\
&=&-g!\ (1+{1\over 2}\ari{c}_1(\mtr{\T\CA}_0))\cdot\Big({(-1)^g\over g!}\deg(Q)+
{(-1)^{g+1}\over (g+1)!}(\ari{c}_1(\mtr{\T\CA}_0)
+[\sum_{l=1}^{g-1}\CH_l])\Big)\\
&=&
-g!\ (1+{1\over 2}\ari{c}_1(\mtr{\T\CA}_0))\cdot\Big({(-1)^g\over g!}+{(-1)^{g+1}\over (g+1)!}(\ari{c}_1(\mtr{\T\CA}_0)
+[\sum_{l=1}^{g-1}\CH_l])\Big)\\
&=&
-(1+{1\over 2}\ari{c}_1(\mtr{\T\CA}_0))\cdot\Big({(-1)^g}+{(-1)^{g+1}\over (g+1)}(\ari{c}_1(\mtr{\T\CA}_0)
+[\sum_{l=1}^{g-1}\CH_l])\Big)\\
\end{eqnarray*}
Si l'on combine ce dernier calcul avec le Lemme \ref{IDKLem}, il vient 
\begin{eqnarray*}
T(\Theta,\mtr{\CO}_\Theta)&=&
-g!\int_{P/B}\Td^{-1}(\mtr{Q})\widetilde{\Td}(\mtr{\cal E})\\
&-&(-1)^g({g+3\over 2g+2})\ari{c}_1(\mtr{\T\CA}_0)-{(-1)^g\over (g+1)}\sum_{l=1}^{g-1}\CH_l
\\
&+& 
g!\int_{P/B}\R(E)\Td(E)-\log({g(g-2)\cdots (1+(1+(-1)^g)/2)\over 
(g-1)(g-3)\cdots(1+(1-(-1)^g)/2)}).
\end{eqnarray*}
Enfin, on a le 

\begin{lemme}
L'égalité 
\begin{align*}
\int_{P/B} \op{Td}^{-1}(\overline{Q})
\widetilde{\op{Td}}(\overline{\C{E}}) 
= (-1)^{g} 
\sum_{p = 0}^{[\frac{g}{2} - 1]}\frac{\zeta_{\mQ}(-1 - 2p)}{(2p + 1)!\,(g - 2p
-1)!}
\C{H}_{2p + 1}
\end{align*}
est vérifiée.
\end{lemme}
Ici $\zeta_{\mQ}(\cdot)$ désigne la fonction z\^eta de Riemann. 
\beginProof
Il suffit de démontrer l'égalité dans le cas où $S=\Spec\ \mC$ et $B$ est un point. 
Nous le supposerons donc pendant la preuve du lemme. 

\`A toute s\'erie formelle sym\'etrique $\phi$ on peut associer 
une classe caract\'eristique encore not\'ee $\phi$ et une classe 
secondaire de Bott-Chern $\tilde{\phi}$ comme dans \cite{Gillet-Soule-CharI}*{\S 1}. 

D'apr\`es le th\'eor\`eme fondamental des fonctions 
sym\'etriques, pour tout fibr\'e vectoriel $F$
la classe $\phi(F)$ s'exprime comme combinaison lin\'eaire finie
des classes de Chern de $F$~:  
\[
\phi(F) = \sum \phi_{i_{1},\dots,i_{k}} c_{i_{1}}(F) \dotsm c_{i_{k}}(F). 
\]
Soit $\overline{\C{F}}$ une suite exacte courte de fibr\'es vectoriels hermitiens
sur une vari\'et\'e complexe $X$~: 
\[
\overline{\C{F}} : 0 \rightarrow \overline{F}' \rightarrow 
\overline{F} \rightarrow \overline{F}'' \rightarrow 0
\]
telle que les m\'etriques sur $F'$ et $F''$ 
sont induites par celle sur $F$. On d\'eduit de 
cette derni\`ere propri\'et\'e que $c(\overline{F}' \oplus \overline{F}'') 
= c(\overline{F})$ en tant que formes, ce qui en appliquant \cite{Gillet-Soule-CharI}*{Proposition 1.3.1} 
\`a chacun des termes monomiaux de~: 
\[
\tilde{\phi}(\overline{\C{F}}) = 
\sum \phi_{i_{1},\dots,i_{k}} \widetilde{c_{i_{1}} \dotsm c_{i_{k}}}(\overline{\C{F}}),
\] 
montre dans $\widetilde{A}(X) := A(X)/ (\op{Im}\partial + 
\op{Im}\overline{\partial})$
l'\'egalit\'e~: 
\[ 
\tilde{\phi}(\overline{\C{F}}) = 
\sum \phi_{i_{1},\dots,i_{k}} \sum_{q = 1}^{k} c_{i_{1}}(\overline{F}) \dotsm 
c_{i_{q-1}}(\overline{F}) \cdot \tilde{c}_{q}(\overline{\C{F}}) \cdot 
c_{i_{q+1}}(\overline{F}) \dotsm
c_{i_{k}}(\overline{F}). 
\]

Appliquons ce qui pr\'ec\`ede pour la classe de Todd et
la suite exacte universelle $\mtr{\cal E}$. Il vient, en remarquant de plus que 
$c(p^*\mtr{\T A}_0) = 1$, l'\'egalit\'e~: 
\[
\widetilde{\op{Td}}(\overline{\C{E}}) = 
\sum_{k = 1}^{g} \op{Td}_{k}
\tilde{c}_{k}(\overline{\C{E}}). 
\]
On sait par ailleurs \cite{Hirzebruch-Topological}*{p.14 Remark 2} que 
$\op{Td}_{k} = \op{Td}_{1,\dots,1}$ ($k$ fois). On peut donc 
\'ecrire~: 
\[
\widetilde{\op{Td}}(\overline{\C{E}}) =
\frac{1}{2}\,\tilde{c}_{1}(\overline{\C{E}}) + 
\sum_{k = 2}^{g} \frac{{B}_{k}}{k!}
\,\tilde{c}_{k}(\overline{\C{E}}),
\]
o\`u $B_{k} = - k \,\zeta_{\M{Q}}(1 - k)$ pour $k \geqslant 2$ est
le $k$-i\`eme nombre de Bernoulli.

On tire de \cite{Gillet-Soule-CharII}*{Proposition 5.3} que 
$\tilde{c}_{1}(\overline{\C{E}}) = 0$ et que pour 
$2 \leqslant k \leqslant g$ on a~:
\[
\tilde{c}_{k}(\overline{\C{E}}) = (- 1)^{k - 1} \C{H}_{k - 1}\, c_{1}^{k - 1}
(\overline{Q}),
\] 
o\`u $\C{H}_{k - 1}$ est le $(k - 1)$-i\`eme nombre harmonique. 

Mettant ce qui pr\'ec\`ede bout-à-bout en tenant compte 
de la nullit\'e des nombres $B_{2p + 1}$ pour $p > 0$, 
on peut finalement \'ecrire~: 
\begin{align*}
\int_{\M{P}^{g - 1}(\M{C})} & \op{Td}^{-1}(\overline{Q})
\widetilde{\op{Td}}(\overline{\C{E}}) \\
= 
&\int_{\M{P}^{g - 1}(\M{C})}
\left(\sum_{q \geqslant 0} \frac{\left(- c_{1}(\overline{Q})\right)^{q}}{(q + 1)!}
\right)
\sum_{p = 0}^{[\frac{g}{2} - 1]}\frac{\zeta_{\M{Q}}(-1 - 2p)}{(2p + 1)!}
\C{H}_{2p + 1}\, c_{1}^{2p + 1}(\overline{Q}) \\
= & (-1)^{g} 
\sum_{p = 0}^{[\frac{g}{2} - 1]}\frac{\zeta_{\M{Q}}(-1 - 2p)}{(2p + 1)!\,(g - 2p
-1)!}
\C{H}_{2p + 1}. 
\end{align*}

\endProof

Enfin, en utilisant la suite exacte universelle $\cal E$, on calcule que 
\begin{eqnarray*}
\int_{P/B}\R(E)\Td(E)&=&-\int_{P/B}\R(Q)\Td^{-1}(Q)=\\
&=&
\int_{P/B}\Big(\big(\sum_{l\geqslant 0}(-1)^{l+1}{x^l\over (l+1)!}\big)\cdot\big(\sum_{m\geqslant 1,\ m\ {\rm impair}}
(2\zeta'_{\mQ}(-m)+\CH_m\zeta_{\mQ}(-m))\cdot{x^m\over m!}\big)\Big)
\end{eqnarray*}
où l'on a posé $x=c_1(Q)$. 
% VERIFIER CE CALCUL
Ainsi
\begin{eqnarray*}
\int_{P/B}\R(E)\Td(E)=
-\sum_{k=0}^{[{g\over 2}-1]}(-1)^{g-2k}{2\zeta'_{\mQ}(-1-2k)+\zeta_{\mQ}(-1-2k)\CH_{2k+1}
\over (2k+1)!(g-2k-1)!}.
\end{eqnarray*}
Le théorème suivant résume maintenant nos calculs:
\begin{theor}
L'égalité
\begin{eqnarray*}
T(\Theta,\mtr{\CO}_\Theta)&=&
(-1)^{g}({g+3\over 2g+2})\ari{c}_1(\mtr{\omega})-{(-1)^g\over (g+1)}\sum_{l=1}^{g-1}\CH_l\\
&-& 
2(-1)^g \sum_{k=0}^{[{g\over 2}-1]}{g\choose 2k+1}\big(\zeta'_{\mQ}(-1-2k)+\zeta_{\mQ}(-1-2k)\CH_{2k+1}\big)\\
&-&\log({g(g-2)\cdots (1+(1+(-1)^g)/2)\over 
(g-1)(g-3)\cdots(1+(1-(-1)^g)/2)})
\end{eqnarray*}
est vérifiée.
\end{theor}
Cette dernière égalité implique notamment  qu'il existe un nombre entier 
$l\in\mN^*$, un nombre réel $C$ et une forme modulaire 
$\mu$ pour le groupe d'Igusa $\Gamma(1,2)$, tels que 
\begin{eqnarray}
\exp(T(\Theta,\mtr{\CO}_\Theta))^l=\Vert C\cdot\mu\Vert^{2l(-1)^{g+1}(g+3)/(2g+2)}_{\rm Pet}
\label{eqY}
\end{eqnarray}
et que $\mu$ est définie sur $\mQ$. 
% pour démontrer cela, il faut utiliser l'astuce de Faltings des espaces 
% de modules avec rigidification linéaire
L'égalité \refeq{eqY} est une forme affaiblie de la première partie de la Conjecture 6.1 de Yoshikawa dans \cite{Yoshikawa-Discriminant}. Une autre 
conséquence de l'égalité \refeq{eqY} est une forme affaiblie  de la première assertion du \og Main Theorem\fg\ dans l'introduction de \cite{Yoshikawa-Discriminant}. 
Dans les deux cas, il s'agit d'une forme affaiblie parce que  le nombre $l$ n'est pas 
effectif et que le groupe d'Igusa $\Gamma(1,2)$ est plus 
petit que le  groupe de Siegel.

{\bf Remarque}. Il est probable que le nombre $l$ peut être déterminé 
en faisant usage dans les calculs ci-dessus du théorème d'Adams-Riemann-Roch arithmétique démontré dans 
\cite{Rossler-Adams} 
(qui tient compte des phénomènes de torsion) plutôt que du théorème de Riemann-Roch arithmétique. Nous n'avons cependant pas effectué ce calcul. 

\section{Les formes modulaires d'Igusa \og produit des thêta constantes paires\fg\  }

Les formes modulaires décrites via la torsion analytique de $\CO_\Theta$ 
dans la dernière section coïncident avec la forme modulaire d'Igusa 
\og produits des fonctions thêta paires\fg\   lorsque $g=2,3$ 
(cf. \cite[après le \og Main Theorem\fg\  ]{Yoshikawa-Discriminant}). 
On peut se demander si ces dernières formes modulaires peuvent aussi 
être interprétées via le théorème de Riemann-Roch arithmétique. 
Nous allons montrer dans cette section qu'une pareille interprétation est 
possible. Le théorème de Riemann-Roch arithmétique est ici appliqué 
à $\CO(\Theta)$ - il s'agit alors d'un cas particulier de la 
\og formule clé\fg\ de Moret-Bailly. 

On continue avec les mêmes hypothèses que dans la section 2. 
 On suppose de plus que  
$\Theta$ est symétrique, i.e. invariant par l'action de l'inversion 
$[-1]$ du schéma en groupes $\CA/B$. Par ailleurs, on suppose que 
$2$ est inversible sur $B$.  

Le théorème du cube (cf. par ex. 
\cite[4.1.23]{Bost-Intrinsic}) implique que
$$
8\cdot\ac1(\mtr{\CO}(\Theta)|_{\CA_{[-1]}})=8\cdot\ac1(u^*\mtr{\CO}(\Theta)).
$$
Par ailleurs, la formule clé de Moret-Bailly (cf. \cite[Th. 3.3]{Moret-Bailly-Sur}) dit que 
$$
8\cdot\ac1(u^*\mtr{\CO}(\Theta))=4\cdot\ac1(\mtr{\omega})+2g\log(4\pi).
$$
On en déduit que $8\cdot\ac1(\mtr{\CO}(\Theta)|_{\CA_{[-1]}})=4\cdot\ac1(\mtr{\omega})$. 
Comme $\Theta$ est lisse sur $B$, le schéma 
$\Theta_{[-1]}$ des points fixes de $[-1]$ dans $\Theta$ est régulier et donc ouvert dans $\CA_{[-1]}$. 
Soit $g:\CA_{[-1]}\backslash\Theta_{[-1]}\to B$ le morphisme de structure. 
% ELUCIDER CET ARGUMENT
On cherche à calculer 
$$
g_*(\ac1(\mtr{\CO}(\Theta)|_{\CA_{[-1]}\backslash\Theta_{[-1]}})).
$$
Vu que le fibré en droites $\CO(\Theta)$ est canoniquement trivialisé sur 
$\CA_{[-1]}\backslash\Theta_{[-1]}$, on est ramené à un calcul 
analytique sur une fibre complexe arbitraire de $\CA\to B$. Nous en fixons 
une et la nommons $A$. 
Nous rappelons 
l'expression explicite de la métrique de Moret-Bailly de $\CO(\Theta)$ donnée dans 
\cite[Par. 3, (3.2.2)]{Moret-Bailly-Sur}. 
Si $\Omega$ est une matrice $g\times g$ complexe du demi-plan de 
Siegel représentant $A$, on dispose de la formule
$$
\Vert s_\Theta(z)\Vert=\det(\Im(\Omega))^{1/4}\exp(-\pi^t y(\Im(\Omega))^{-1}y)|\theta(z,\Omega)|
$$
où $\theta(z,\Omega)$ est la fonction $\theta$ de Riemann associée à $\Omega$, 
$z=x+iy$ et $s_\Theta$ est la section canonique de $\CO(\Theta)$, restreinte 
à $A$. Si on utilise la formule de changement de coordonnées
$$
x+iy=x_1+\Omega(y_1)=(x-\Re(\Omega)(\Im(\Omega))^{-1}y)+\Omega((\Im(\Omega))^{-1}y)
$$
on peut réécrire, en utilisant la symétrie de $\Im(\Omega)$, 
$$
\Vert s_\Theta(z)\Vert=\det(\Im(\Omega))^{1/4}\exp(-\pi^t y_1\Im(\Omega)y_1)|\theta(z,\Omega)|.
$$
On calcule maintenant 
\begin{eqnarray*}
g_*(\ac1(\mtr{\CO}(\Theta)|_{\CA_{[-1]}\backslash\Theta_{[-1]}}))&=&
-2\log|\prod_{a,b}\det(\Im(\Omega))^{1/4}\exp(-\pi^t a\Im(\Omega)a)\theta(
\Omega(a)+b,\Omega)|\\
&=&-\log|\det(\Im(\Omega))|^{2^{2g-2}+2^{g-2}}-\log|\prod\theta\left[{a\atop b}\right](0,\Omega)|^2\\
&=:&
-\log\Vert\chi_g(\Omega)\Vert^2_{\rm Pet},
\end{eqnarray*}
où $\chi_g(\cdot)$ est la forme modulaire d'Igusa 
\og produit des thêta constantes paires\fg\ 
(cf. \cite[chap. II]{Mumford-Tata-I}) et $\Vert\cdot\Vert_\Pet$ est la norme associée 
à la métrique de Petersson. 
La somme porte sur les couples $(a,b)\in\mQ^g$ tels 
que $a,b\in\{0,1/2\}$ et tels que $\theta\left[{a\atop b}\right](0,\Omega)\not= 0$. 
Il y a $2^{2g}-2^{g-1}(2^g-1)=2^{2g-1}+2^{g-1}$ tels couples. Ceci résulte par exemple 
de la formule de Lefschetz habituelle appliquée à $[-1]$ agissant sur $\Theta_\mC|_A$. 
Par ailleurs,
$$
8\cdot g_*(\ac1(\mtr{\CO}(\Theta)|_{\CA_{[-1]}\backslash\Theta_{[-1]}}))
=
4\cdot(2^{2g-1}+2^{g-1})\ac1(\mtr{\omega})+2g\cdot(2^{2g-1}+2^{g-1})\log(4\pi)
$$
ce qui implique la 
\begin{prop}
L'égalité
$$
(2^{2g-2}+2^{g-2})\ac1(\mtr{\omega})+{g\over 4}(2^{2g-1}+2^{g-1})\log(4\pi)=-\log\Vert\chi_g(\Omega)\Vert^2_{\rm Pet}
$$
est vérifiée.
\end{prop}
Supposons \`a pr\'esent que
$S$ est le spectre d'un anneau de Dedekind et affaiblissons
les hypoth\`eses pr\'ec\'edentes en supposant seulement 
que $\Theta$ est lisse au-dessus d'un ouvert non-vide $V$ de $B$. 
Supposons de plus que $B=S$  et 
que le schéma $\CA_\iota$ est la réunion disjointe de 
$2^{2g}$ sections $B\to\CA$.  

Soit maintenant 
$$
Z:=\Zar(\CA_{[-1],V}\backslash\Theta_{[-1],V})
$$
l'adhérence schématique de $\CA_{[-1],V}\backslash\Theta_{[-1],V}$. 
On 
cherche alors à calculer
$$
g_*(\ac1(\mtr{\CO}(\Theta)|_Z)).
$$
Soit $u_1,\dots,u_{(2^{2g-1}+2^{g-1})}$ une énumération des sections formant 
$Z$. 
Une variante  du calcul fait plus plus haut donne alors la
\begin{prop}
L'égalité 
\begin{eqnarray*}
(2^{2g-2}+2^{g-2})\ac1(\mtr{\omega})&+&{g\over 4}(2^{2g-1}+2^{g-1})\log(4\pi)\\
&=&\sum_{j=1}^{2^{2g-1}+2^{g-1}}\sum_{P\in 
u_j\cap\Theta}{\rm long}_{\CO_{\CA,P}}(\CO_{u_j\cap\Theta})\cdot\pi_*(P)-\log\Vert\chi_g(\Omega)\Vert^2_{\rm Pet}
\end{eqnarray*}
est vérifiée.
\end{prop}

\section{Les formes modulaires de Yoshikawa de deuxième type}

On se donne un schéma régulier et intègre $\CB$, quasi-projectif 
sur un anneau arithmétique 
$D$ de corps de fractions $K$. 
On se donne également un schéma intègre et régulier $\CT$ et 
un morphisme projectif, plat $f:\CT\ra \CB$ tel que $f_K:\CT_K\ra \CB_K$ est lisse.
On notera $B:=\CB_K$ (resp. $T:=\CT_K$) la fibre générique de 
$B$ (resp. $\CT$) sur $D$. 

On suppose que $T$ est un schéma en surfaces K3 sur $B$. 
Par définition, cela signifie que les fibres géométriques de $f_K$ sont des surfaces 
K3. On munit $T_\mC$ d'une structure de fibration Kählerienne $\nu$. 
On suppose aussi  que le morphisme d'adjonction
$$
f^* f_{*}{\omega}\ra\omega
$$
est un isomorphisme.  Munissons $\omega$ de la métrique induite par 
la structure de fibration Kählerienne et $f^* f_{*}{\omega}$ de la métrique image réciproque par $f$ de la métrique $L^2$. 
On écrira $\eta$ pour la classe secondaire $\widetilde{\ch}(\mtr{\cal F})$ 
de la suite exacte de fibrés
$$
{\cal F}:0\ra f^* f_{*}{\omega}\ra\omega\ra 0\ra 0
$$
munie de ces métriques. 
On notera $\eta_g$ la restriction de cette classe à $\CT_{g,\mC}$. 
On suppose aussi qu'il existe un 
automorphisme d'ordre $2$ de $\CT$ sur $\CB$. Par ailleurs, on 
suppose que $2$ est inversible sur $D$. On dispose ainsi d'une action 
de $\mu_2$ sur $\CT$. On note $g$ l'action de l'automorphisme 
sur $T(\mC)$. 
On suppose également que la fibration Kählerienne est équivariante 
pour l'action de $g$ et que le morphisme 
$\CT_\m2\ra\CB$ est lisse. 

On applique la formule de Lefschetz arithmétique au morphisme $f$, à l'action de $\m2$ sur $\CT$ et au fibré trivial $\CO:=\CO_\CT$ muni de sa métrique triviale.
 On note $N$ le fibré conormal de 
$\CT_{\m2}$ dans $\CT$ et on le munit de la métrique 
induite par la structure de fibration Kählerienne. 
On obtient 
\begin{eqnarray*}
&&\ari{c}_{1,\m2}(R^0 f_*\mtr{\CO})-\ari{c}_{1,\m2}(R^1 f_*\mtr{\CO})+\ari{c}_{1,\m2}(R^2 f_*\mtr{\CO})\\
&=& f_*(\ari{\Td}_\m2(\mtr{Tf}))^{[1]}-
\int_{T_\m2/B}\Td_g(Tf_\mC)R_g(Tf_\mC)+T_g(\mtr{\CO}).
\end{eqnarray*}
On calcule dans $\ari{\rm CH}^{\leqslant 2}(\CT_\m2)_\mQ$:
\begin{eqnarray*}
\ari{\Td}_\m2(\mtr{Tf})&=&\ari{\ch}_\m2(1-\mtr{N})^{-1}\ari{\Td}(\mtr{Tf}_\m2)\\
&=&
{1\over 2}(1+{1\over 2}\ac1(\mtr{N})+{1\over 4}\ari{c}_1(\mtr{N})^2)^{-1}\ari{\Td}(\mtr{Tf}_\m2)\\
&=&
\Big({1\over 2}-{1\over 4}\ac1(\mtr{N})\Big)\ari{\Td}(\mtr{Tf}_\m2)
\end{eqnarray*}
Par ailleurs, dans $\ari{\rm CH}^{\leqslant 2}(\CT_\m2)_\mQ$, on a 
 $\ari{\Td}(\mtr{Tf}_\m2)=1-{1\over 2}\ac1(\mtr{\omega}_\m2)+{1\over 12}\ac1(\mtr{\omega}_\m2)^2$, 
où $\mtr{\omega}_\m2$ est le fibré des différentielles relative  
de $\CT_\m2$ sur $\CB$, muni de la métrique induite. 
La partie de degré 
$2$ de $\ari{\Td}_\m2(\mtr{Tf})$ est donc la partie de 
degré $2$ de l'expression
$$
\Big({1\over 2}-{1\over 4}\ac1(\mtr{N})\Big)\Big(1-{1\over 2}\ac1(\mtr{\omega}_\m2)
+{1\over 12}\ac1(\mtr{\omega}_\m2)^2\Big)
$$
qui est 
\begin{equation}
{1\over 8}\ac1(\mtr{N})\ac1(\mtr{\omega}_\m2)+{1\over 24}\ac1(\mtr{\omega}_\m2)^2.
\label{SS}
\end{equation}
On rappelle qu'on dispose d'une suite exacte équivariante
$$
0\ra N\ra\Omega\ra\omega_\m2\ra 0
$$
sur $\CT_\m2$. Pour des raisons de rang, cette suite est isométriquement 
scindée. On a donc 
$$
\ac1(N)=f^*_\m2\ac1(f_*\mtr{\omega})-\eta_g-\ac1(\mtr{\omega}_\m2)
$$
dans $\ari{\rm CH}^1(\CT_\m2)$. On peut donc évaluer l'expression \refeq{SS} 
comme
$$
{1\over 8}\Big(f^*_\m2\ac1(f_*\mtr{\omega})-\eta_g-\ac1(\mtr{\omega}_\m2)\Big)\ac1(\mtr{\omega}_g)
+{1\over 24}\ac1(\mtr{\omega}_\m2)^2
={1\over 8}f^*_\m2\ac1(f_*\mtr{\omega})\ac1(\mtr{\omega}_\m2)-{1\over 12}\ac1(\mtr{\omega}_\m2)^2-{1\over 8}c_1(\mtr{\omega}_\m2)\eta_g.
$$
Par ailleurs, on calcule
\begin{eqnarray*}
&&\int_{T_g/B}\Td_g(Tf_\mC)R_g(Tf_\mC)\\
&=&-2\int_{T_g/B}\Big((2\zeta'_{\mQ}(-1,-1)+\zeta_{\mQ}(-1,-1))c_1(N)+(2\zeta'_{\mQ}(-1)+\zeta_{\mQ}(-1))c_1(\omega_g)\Big)\\
&=&-2\int_{T_g/B}\Big((6\zeta'_{\mQ}(-1)+(3-\log(16))\zeta_{\mQ}(-1))c_1(N)+(2\zeta'_{\mQ}(-1)+\zeta_{\mQ}(-1))c_1(\omega_g)\Big)\\
&=&-2\int_{T_g/B}\Big(-(6\zeta'_{\mQ}(-1)+(3-\log(16))\zeta_{\mQ}(-1))+(2\zeta'_{\mQ}(-1)+\zeta_{\mQ}(-1))\Big)c_1(\omega_g)\\
&=&
-2\int_{T_g/B}\Big(-4\zeta'_{\mQ}(-1)+(\log(16)-2)\zeta_{\mQ}(-1)\Big)c_1(\omega_g)\\
&=&-2G\Big(-4\zeta'_{\mQ}(-1)+(\log(16)-2)\zeta_{\mQ}(-1)\Big)
\end{eqnarray*}
où $G$ est une fonction localement constante sur $B(\mC)$. En un point 
$P\in B(\mC)$, $G$ vaut 
$$
\sum_{C\subseteq T_{g,P}(\mC)}(2\cdot{\rm genre}(C)-2)
$$
où la somme porte sur les composantes connexes $C$ de la fibre 
$T_{g,P}(\mC)$ de $T_{g}(\mC)$ au-dessus de $P$. 
Pour résumer, on obtient
\begin{eqnarray*}
&&\ari{c}_{1,\m2}(R^0 f_*\mtr{\CO})-\ari{c}_{1,\m2}(R^1 f_*\mtr{\CO})+\ari{c}_{1,\m2}(R^2 f_*\mtr{\CO})\\
&=&
{G\over 8}\ac1(f_*\mtr{\omega})-{1\over 12}f_{\m2*}\ac1(\mtr{\omega}_\m2)^2
+2G\Big(-4\zeta'_{\mQ}(-1)+(\log(16)-2)\zeta_{\mQ}(-1)\Big)\\
&+&T_g(\mtr{\cal O})-{1\over 8}\int_{T_g/B}
c_1(\mtr{\omega}_g)\eta_g.
\end{eqnarray*}
Ceci implique en particulier le 
\begin{theor}
Supposons que toutes les fibres géométriques de $f$ sont des surfaces K3, 
alors on a 
\begin{eqnarray*}
-\log|{1\over d!({2\pi})^d}\int_{T/B}\nu^d|&=&
{G-8\over 8}\ac1(f_*\mtr{\omega})-{1\over 12}f_*\ac1(\mtr{\omega}_\m2)^2
+2G\Big(-4\zeta'_{\mQ}(-1)+(\log(16)-2)\zeta_{\mQ}(-1)\Big)\\
&+&T_g(\mtr{\cal O})-{1\over 8}\int_{T_g/B}
c_1(\mtr{\omega}_\m2)\eta_g.
\end{eqnarray*}
\label{thcor}
\end{theor}
On remarquera que si $\CT_g$ est vide, on a 
\begin{eqnarray*}
-\log|{1\over d!({2\pi})^d}\int_{\CT/\CB}\nu^d|+\ac1(f_*\mtr{\omega})=T_g(\mtr{\CO}).
\label{enri}
\end{eqnarray*}
sous les hypothèses du théorème \ref{thcor}. 
On peut exprimer la quantité ${1\over 12}f_{g*}\ac1(\mtr{\omega}_\m2)^2$ du Théorème \ref{thcor} 
au moyen de la torsion analytique des fibres de $T_g(\mC)$ sur $B(\mC)$, via le théorème de 
Riemann-Roch arithmétique. On obtient 
\begin{eqnarray*}
{1\over 12}f_{\m2*}\ac1(\mtr{\omega}_\m2)^2=f_*(\ari{\Td}(\CT/\CB))^{[1]}&=&
-T(\mtr{\CO}_{g})-\log|{1\over d_g!(2\pi)^{d_g}}\int_{\CT_g/\CB}\nu_g^{d_g}|+\ac1(f_{\m2*}\mtr{\omega}_\m2)\\&-&\int_{\CT_g/\CB}(2\zeta'_{\mQ}(-1)+\zeta_{\mQ}(-1))c_1(\omega_\m2)
\end{eqnarray*}
dans $\ari{\rm CH}^1(\CB)_\mQ$. Si l'on juxtapose cette dernière expression à celle 
du Théorème \ref{thcor}, on obtient
\begin{eqnarray*}
&&\ac1(f_{g*}\mtr{\omega}_\m2)+{8-G\over 8}\ac1(f_*\mtr{\omega})\\
&=&T_g(\mtr{\cal O})+T(\mtr{\cal O}_g)
-{1\over 8}\int_{T_g/B}
c_1(\mtr{\omega}_\m2)\eta_g+\log|{1\over d!({2\pi})^d}\int_{T/B}\nu^d|+
\log|{1\over d_g!(2\pi)^{d_g}}\int_{T_g/B}\nu_g^{d_g}|\\
&+&2G\Big(-4\zeta'_{\mQ}(-1)+(\log(16)-2)\zeta_{\mQ}(-1)\Big)
+\int_{T_g/B}(2\zeta'_{\mQ}(-1)+\zeta_{\mQ}(-1))c_1(\omega_\m2)\\
&=&T_g(\mtr{\cal O})+T(\mtr{\cal O}_g)
-{1\over 8}\int_{T_g/B}
c_1(\mtr{\omega}_\m2)\eta_g+\log|{1\over d!({2\pi})^d}\int_{T/B}\nu^d|+
\log|{1\over d_g!(2\pi)^{d_g}}\int_{T_g/B}\nu_g^{d_g}|\\
&-&6G\zeta'_{\mQ}(-1)-{2G\over 3}\log(2)+{G\over 4}
\end{eqnarray*}
sous les hypothèses du Théorème \ref{thcor}. 
On suppose maintenant que $D=\mC$; les hypothèses du Théorème \ref{thcor} sont alors 
automatiquement satisfaites. Supposons par ailleurs que 
$f_*\omega$ a une section analytique trivialisante de norme $L^2$ constante. 
Ceci est le cas par exemple si la famille $\CT$ est munie d'un marquage 
(cf. \cite[Par. 1.2 (b)]{Yoshikawa-K3-surfaces} pour cette notion). 
% (on se rappellera du théorème de Griffiths disant qu'une section localement 
% constante d'une variation de structure de Hodge de type pure en un point 
% l'est en tout point [cf. ... Deligne - Hodge II]). 
Soit $\kappa$ un entier tel que le fibré
 $f_*\omega^{\otimes (8-G)\kappa}\otimes
 (\det f_{\m2*}\omega_\m2)^{\otimes 8\kappa}$ est trivial. 
Le fibré $(\det f_{\m2*}\omega_\m2)^{\otimes(-8\kappa)}$ est alors 
analytiquement trivial. Il existe donc $t$ une section analytique trivialisante 
de $(\det f_{g*}\omega_g)^{\otimes 8\kappa}$ satisfaisant l'égalité
\begin{eqnarray*}
|t|_{L^2}^{-{1\over 4\kappa}}=e^{T_g(\mtr{\cal O})}\cdot e^{T(\mtr{O}_g)}\cdot
|{1\over d!({2\pi})^d}\int_{T/B}\nu^d|\cdot
|{1\over d_g!(2\pi)^{d_g}}\int_{T_g/B}\nu_g^{d_g}|\cdot
\exp(-{1\over 8}\int_{T_g/B}
c_1(\mtr{\omega}_\m2)\eta_g)
\end{eqnarray*}
\'Ecrivons $\Vol(T_g):=|{1\over d_g!(2\pi)^{d_g}}\int_{T_g/B}\nu_g^{d_g}|$ 
et $\Vol(T):=|{1\over d!({2\pi})^d}\int_{T/B}\nu^d|$. Soit 
$r_+$ (resp. $r_{-}$ la dimension du sous-espace de $H^2(T(\mC)_b,\mC)$ invariant par 
$\iota$ (resp. celui où $\iota$ agit par $-1$); $b$ est un élément générique 
de $B(\mC)$. 
% Dans la notation de Yoshikawa, on $r(M)=r_+$. 
Remarquons que 
par la formule du point fixe holomorphe et la formule de Gauss-Bonnet généralisée 
(cf. \cite[Example 3.8, chap. III, sec. 3, p.96]{Wells-Differential} pour cette dernière), on a l'égalité
$$
1-0+r_{+}-r_{-}+0-1=-G
$$
et par ailleurs, le formulaire \cite[VIII, 3.]{Barth-Hulek-Compact} nous assure que 
$r_{+}+r_{-}=22$. 
On en déduit que
$$
G=20-2r_+.
$$
On reprend maintenant l'expression pour $|t|_{L^2}^{-{1\over 4\kappa}}$ et 
on calcule 
\begin{eqnarray*}
&&e^{T_g(\mtr{\cal O})}\cdot e^{T(\mtr{O}_g)}\cdot
\Vol(T)\cdot\Vol(T_g)\cdot
\exp(-{1\over 8}\int_{\CT_g/\CB}
c_1(\mtr{\omega}_\m2)\eta_g)\\
&=&
e^{T_g(\mtr{\cal O})}\cdot e^{T(\mtr{O}_g)}\cdot
\Vol(T)\cdot\Vol(T_g)\cdot
\exp\big(-{1\over 8}\int_{T_g/B}
c_1(\mtr{\omega}_\m2)(\eta_g+\log|\Vol(T)|)\big)\cdot\Vol(T)^{G\over 8}\\
&=&
e^{T_g(\mtr{\cal O})}\cdot e^{T(\mtr{O}_g)}\cdot
\Vol(T)^{G/8+1}\cdot\Vol(T_g)\cdot
\exp\big(-{1\over 8}\int_{T_g/B}
c_1(\mtr{\omega}_\m2)(\eta_g+\log|\Vol(T)|)\big).
\end{eqnarray*}
Par ailleurs, on a 
$$
G/8+1={G+8\over 8}={20-2r_+ +8\over 8}={14-r_+\over 4}
$$
et on conclut que 
\begin{eqnarray*}
|t|_{L^2}^{-{1\over 4\kappa}}=
e^{T_g(\mtr{\cal O})}\cdot e^{T(\mtr{O}_g)}\cdot
\Vol(T)^{14-r_+\over 4}\cdot\Vol(T_g)\cdot
\exp\big(-{1\over 8}\int_{T_g/B}
c_1(\mtr{\omega}_\m2)(\eta_g+\log|\Vol(T)|)\big)
\end{eqnarray*}
Il s'agit  
de l'égalité du théorème principal \cite[Main Th., Introduction]{Yoshikawa-K3-surfaces} 
de Yoshikawa.

\section{Appendice: une formule du point fixe singulière conjecturale en théorie 
d'Arakelov}

Soit $D$ un anneau arithmétique d'anneau de fractions $K$ 
et supposons que $D$ est régulier. 
Soit $f:X\to\Spec\ D$ un schéma intègre, projectif sur $D$, dont la fibre sur 
$K$ est lisse. 
Soit $h:Z\to\Spec\ D$ un schéma intègre et régulier, projectif sur $D$, dont la fibre sur 
$K$ est lisse. 
Soit $j$ une $D$-immersion fermée $X\hookrightarrow Z$. 
On munit $Z$ d'une métrique Kählerienne $\omega_Z$ et 
on munit $X$ de la structure $\omega_X$ induite. 
On se donne un nombre entier $n\geqslant 1$ et 
des structures $\mn$-équivariantes sur $X$ et $Z$ telle que 
$f,h,j$ soient $\mn$-équivariants et que la structure $\omega_Z$ soit 
$\mn(\mC)$-invariante ($D$ est supposé muni de la structure 
équivariante triviale). 
Soit enfin $R(\mn)=\mZ/(1-T^n)$ le groupe de Grothendieck 
des $\mn$-comodules de type fini sur $\mZ$. 
On choisit un racine primitive $n$-ième 
de l'unité $\zeta_n$ et 
 une $R(\mn)$-algèbre $\cal R$ telle que les éléments  
$1-T^{k}$ ($k=1,\dots,n-1$) sont inversibles dans $\cal R$.

Soit $N$ le fibré conormal de l'immersion $Z_\mn\hookrightarrow Z$. 
Soit enfin $\sum_{i}r_i\mtr{E}_i$ une $\cal R$ combinaison linéaire 
finie de fibrés hermitiens sur $Z_\mn$ tels que 
$\sum_i{r_i}\mtr{E}_i=(\Lambda_{-1}(\mtr{N}))^{-1}$ dans 
$\widehat{\K}_0^\mn(Z)\otimes_{\R(\mn)}{\cal R}$. 

Soit $\mtr{E}$ un fibré hermitien $\mn$-équivariant sur $X$. 

On remarque que l'immersion $X_\mC\hookrightarrow Z_\mC$ est régulière 
et on a donc 
$$\ul{\Tor}^k_{\CO_Z}(j_*E,\CO_{Z_\mn})_\mC\simeq 
j_{\mC*}(\Lambda^k(F)\otimes E_\mC),
$$
 où $F$ est un fibré localement 
libre défini sur $X_{\mn,\mC}$ par la suite exacte
$$
\CF: 0\to F\to N_{Z_{\mn,\mC}/Z_{\mC}}\to N_{X_{\mn,\mC}/X_{\mC}}\to 0
$$
(voir \cite[Exp. VII, Prop. 2.5]{Grothendieck-SGA6}). Nous munissons le fibré $F$ de la métrique induite par $N_{Z_\mn/Z}$.

Pour tout $l\geqslant 0$, les fibrés cohérents $\R^l h_*(E_i\otimes \ul{\Tor}^k_{\CO_Z}(j_*E,\CO_{Z_\mn}))$ (qui sont localement libres sur la fibre générique) 
peuvent être munis de métriques hermitiennes via 
l'isomorphisme naturel
$$
\R^l h_*(E_i\otimes \ul{\Tor}^k_{\CO_Z}(j_*E,\CO_{Z_\mn}))_\mC
\simeq \R^l f_{\mC*}(j^*(E_{i,\mC})\otimes\Lambda^k(F)\otimes E_{\mC})
$$
Par abus de notation, on notera $\R^l h_*(\mtr{E}_i\otimes\ul{\Tor}^k_{{\CO}_Z}(j_*\mtr{E},\mtr{\CO}_{Z_\mn}))$ le fibré cohérent hermitien sur $D$ 
(\og hermitian coherent sheaf\fg\   en anglais) correspondant. 

\begin{conj}
L'égalité 
\begin{eqnarray*}
\sum_{l\geqslant 0}(-1)^l\R^l f_*(\mtr{E})-\T_g(\mtr{E})&=&\sum_{i}r_i\sum_{l,k\geqslant 0}(-1)^{l+k}\R^l h_*(\mtr{E}_i\otimes \ul{\Tor}^k_{\CO_Z}(j_*\mtr{E},\mtr{\CO}_{Z_\mn}))\\
&+&\int_{X_{\mn}}\ch_g(\mtr{E}_\mC)\widetilde{\Td}_g(\mtr{\CF})\Td^{-1}_g(\mtr{F})
\\&-&\int_{X_\mn}\ch_g({E}_\mC)\Td_g(\T X_{\mC})\R_g(N_{X_{\mn,\mC}/X_{\mC}})
\end{eqnarray*}
est vérifiée dans $\widehat{\K}^{\mn'}_0(D)\otimes_{\R(\mn)}{\cal R}$.
\end{conj}

Cette conjecture est inspirée par la formule \cite[Th. 3.5]{Thomason-Lefschetz}.

\end{document}